# Extended Conway-Maxwell-Poisson distribution and its properties and applications

*Subrata Chakraborty[1], Tomoaki Imoto[2]*

*[1]Department of Statistics, Dibrugarh University, Dibrugarh-786004, Assam, India*

*Email: subrata_arya@yahoo.co.in*

*[2]The institute of Statistical Mathematics, 10-3 Midori-cho, Tachikawa, Tokyo-190-8562, Japan*

**Abstract**
A new three parameter natural extension of the Conway-Maxwell-Poisson (COM-Poisson) distribution is proposed. This distribution includes the recently proposed COM-Poisson type negative binomial (COM-NB) distribution [*Chakraborty, S. and Ong, S. H. (2014): A COM-type Generalization of the Negative Binomial Distribution, Accepted in Communications in Statistics-Theory and Methods*] and the generalized COM-Poisson (GCOMP) distribution [*Imoto, T. :(2014) A generalized Conway-Maxwell-Poisson distribution which includes the negative binomial distribution, Applied Mathematics and Computation, 247, 824-834*]. The proposed distribution is derived from a queuing system with state dependent arrival and service rates and also from an exponential combination of negative binomial and COM-Poisson distribution. Some distributional, reliability and stochastic ordering properties are investigated. Computational asymptotic approximations, different characterizations, parameter estimation and data fitting example also discussed.

## 1. Introduction

Recently two new generalizations of the well known COM-Poisson (Conway and Maxwell, 1962) was proposed. One by Chakraborty and Ong (2014) known as the COM-Negative binomial distribution and the second one by Imoto (2014) referred to as the generalized COM-Poisson Distribution.

### 1.1 COM-Poisson type negative binomial distribution

Chakraborty and Ong (2014) proposed a new COM-Poisson type generalization of negative binomial distribution that includes some well-known distributions including COM-Poisson (Conway and Maxwell, 1962), Negative Binomial (Johnson et al., 2005), as particular case and Bernoulli (Johnson et al., 2005), COM-Poisson (Conway and Maxwell, 1962) as limiting cases among others. This distribution is log-concave and flexible enough to model under, equi- and over dispersed count data.

A random variable $X$ is said to follow the COM - Poisson type Negative Binomial distribution with parameters $(\nu, p, \alpha)$ [COM-NB$(\nu, p, \alpha)$] if its pmf is given by

$$P(X = k) = (\nu)_k \, p^k / \{(k!)^\alpha \, {}_1H_{\alpha-1}(\nu; 1; p)\}, \, k = 0, 1, 2, \cdots \qquad (1)$$

Where
$$ {}_1H_{\alpha-1}(\nu; 1; p) = \sum_{k=0}^{\infty} (\nu)_k \, p^k / (k!)^\alpha \qquad (2)$$

The distribution is defined in the parameter space
$\{\nu > 0, p > 0, \alpha > 1\} \bigcup \{\nu > 0, 0 < p < 1, \alpha = 1\}$. When $\alpha$ is a positive integer, ${}_1H_{\alpha-1}(\nu; 1; p)$ can be expressed as a particular generalized hypergeometric series

$$ {}_mF_n(a_1, a_2, \cdots, a_m; b_1, b_2, \cdots, b_m; z) = \sum_{k=0}^{\infty} \frac{(a_1)_k (a_2)_k \cdots (a_m)_k}{(b_1)_k (b_2)_k \cdots (b_n)_k} \frac{z^k}{k!} \text{ as } {}_1F_{\alpha-1}(\nu; 1, 1, \cdots, 1; p). $$

*1.2 Generalized COM-Poisson distribution*

Imoto (2014) proposed another generalization where a random variable $X$ is said to follow the GCOM-Poisson distribution of Imoto (2014) with parameters $(\nu, p, \beta)$ that is GCOMP $(\nu, p, \beta)$ if its pmf is given by

$$P(X = k) = \frac{\{\Gamma(\nu + k)\}^\beta}{k!} p^k \Big/ C(\beta, \nu, p) \qquad (3)$$

Where
$$C(\beta, \nu, p) = \sum_{k=0}^{\infty} \frac{\{\Gamma(\nu + k)\}^\beta}{k!} p^k \qquad (4)$$

The distribution is defined in the parameter space
$$\{\nu > 0, p > 0, \beta < 1\} \bigcup \{\nu > 0, 0 < p < 1, \beta = 1\}.$$

In the present article a natural three parameter generalization of the COM-Poisson distribution which includes COM-NB and GCOM-Poisson distributions is proposed. Some important distributional properties of the proposed distribution are presented in the section 2. Reliability and stochastic ordering results are discussed in section3. Concluding remarks is made in the final section.

*1.3 A hypergeometric type series*
We introduce the series

$$ {}_mS_\alpha^\beta(a_1, a_2, \cdots, a_m; b; p) = \sum_{k=0}^{\infty} \frac{\{(a_1)_k\}^\beta (a_2)_k \cdots (a_m)_k}{\{(b)_k\}^\alpha} \frac{p^k}{k!}, $$

where $(a)_k = a(a+1) \cdots (a+k-1) = \Gamma(a+k)/\Gamma a$ is the Pochhammer's notation (see Johnson et al., 2005). For $\alpha, \beta$ and $m$ all positive integers, it reduces to ${}_{\beta+m-1}F_\alpha(a_1, a_1, \cdots, a_1, a_2, \cdots, a_m; b, b, \cdots, b; p)$, a particular generalized hypergeometric function. With this notation we have

$$\sum_{k=0}^{\infty} \{(\nu)_k\}^\beta \, p^k / (k!)^\alpha = {}_1S_{\alpha-1}^\beta(\nu; 1; p) = {}_\beta F_{\alpha-1}(\nu; 1, 1, \cdots, 1; p) \qquad (5)$$

Some important special cases of ${}_1S_{\alpha-1}^\beta(\nu; 1; p)$ are

i.   $_1S_{\alpha-1}^1(\nu;1;p) =_1H_{\alpha-1}(\nu;1;p)$ [Chakraborty and Ong, 2014]

ii.   $_1S_o^\beta(\nu;1;p) = C(\beta,\nu,p)/(\Gamma\nu)^\beta$ [Imoto, 2014]

iii.   $_1S_{\alpha-1}^\beta(1;1;p) = Z(p,\alpha-\beta)$ [Conway and Maxwell, 1962]

iv.   $_1S_0^1(\nu;1;p) = (1-p)^{-\nu}$ [geometric series]

Some important limiting cases of $_1S_{\alpha-1}^\beta(\nu;1;p)$ are

v.   $\lim\limits_{\alpha \to \infty}{}_1S_{\alpha-1}^\beta(\nu;1;p) = 1+\nu^\beta p$ .

vi.   $\lim\limits_{\nu \to \infty}{}_1S_{\alpha-1}^\beta(\nu;1;p) = \sum\limits_{k=0}^{\infty}\lambda^k/(k!)^\alpha = Z(\lambda,\alpha)$ , where $\nu^\beta p = \lambda$ is finite positive.

## 2. Extended COM-Poisson (ECOMP) distribution

Here we introduce a new distribution that unifies both the COM-NB and GCOMP distributions.

**Definition1**. A random variable $X$ is said to follow the extended COM-Poisson distribution with parameters $(\nu,p,\alpha,\beta)$ [ECOMP $(\nu,p,\alpha,\beta)$ ] if its pmf is given by

$$P(X=k) = \frac{\{(\nu)_k\}^\beta}{(k!)^\alpha}p^k \bigg/ {}_1S_{\alpha-1}^\beta(\nu;1;p) \qquad (6)$$

The distribution is defined in the parameter space
$\{\nu \geq 0, p > 0, \alpha > \beta\} \cup \{\nu > 0, 0 < p < 1, \alpha = \beta\}$ .

**Remark1.** Unlike COM-NB where the parameter $\alpha \geq 1$ and in GCOMP the parameter $\beta \leq 1$ , in case of ECOMP these two parameters can be either positive or negative with the restriction of $\alpha \geq \beta$ .

**Remark2.** The pmf in (6) can alternatively expressed as

$$P(X=k) = \frac{\{\Gamma(\nu+k)\}^\beta}{(\Gamma\nu)^\beta(k!)^\alpha}p^k \bigg/ {}_1S_{\alpha-1}^\beta(\nu;1;p)$$

### 2.1 Particular cases

Following discrete distributions are particular cases of ECOMP $(\nu,p,\alpha,\beta)$ :

- ➢ COM-NB $(\nu,p,\alpha)$ for $\beta = 1$
- ➢ GCOMP $(\nu,p,\beta)$ for $\alpha = 1$
- ➢ NB $(\nu,p)$ for $\alpha = \beta = 1$ with pmf $P_n = \binom{\nu+n-1}{n}p^n(1-p)^\nu$
- ➢ COMP $(p,\alpha-\beta)$ for $\nu = 1$
- ➢ COMP $(p,\alpha)$ for $\beta = 0$
- ➢ Poisson $(p)$ for $\nu = 1, \alpha = \beta$
- ➢ Poisson $(p)$ for $\beta = 0, \alpha = 1$

➢ A new generalization of NB(NGNB) distribution when $\alpha = \beta = \gamma$ with pmf

$$P(X = k) = \binom{\nu + k - 1}{k}^{\gamma} p^k \Bigg/ {}_1S^{\gamma}_{\gamma-1}(\nu;1;p) \qquad (7)$$

The distribution in (7) is log-convex as will be seen in section 2.8.

*2.2 Approximation of the normalizing constant*

The normalizing constant ${}_1S^{\beta}_{\alpha-1}(\nu;1;p)$ of the ECOMP $(\nu, p, \alpha, \beta)$ distribution is not expressed in a closed form and includes the summation of infinite series. Therefore, we need approximations of this constant to compute the pmf and moments of the distribution numerically.

A simply approximation is to truncate the series that is

$$ {}_1S^{\beta}_{\alpha-1,m}(\nu;1;p) = \sum_{k=0}^{m} \frac{\{(\nu)_k\}^{\beta}}{(k!)^{\alpha}} p^k, \qquad (8)$$

where $m$ is an integer chosen such that $\varepsilon_m = (\nu - m + 1)^{\beta} p / m^{\alpha} < 1$. The relative truncation error is then given by the expression $R_m(\nu, p, \alpha, \beta)$ $= \{{}_1S^{\beta}_{\alpha-1}(\nu;1;p) - {}_1S^{\beta}_{\alpha-1,m}(\nu;1;p)\} / {}_1S^{\beta}_{\alpha-1,m}(\nu;1;p)$. Then the relative error about the pmf is give by $\{P_m(k) - P(k)\} / P(k)$, where $P(k)$ is given by the r.h.s. of equation (6) in section 2 and $P_m(k)$ is given by the r.h.s. of (6) with ${}_1S^{\beta}_{\alpha-1}(\nu;1;p)$ substituted by ${}_1S^{\beta}_{\alpha-1,m}(\nu;1;p)$. The upper bound of the relative truncation error is then found to be

$$R_m(\nu, p, \alpha, \beta) < \frac{\{(\nu)_{m+1}\}^{\beta} p^{m+1}}{\{(m+1)!\}^{\alpha} {}_1S^{\beta}_{\alpha-1,m}(\nu;1;p)} \sum_{k=0}^{\infty} \varepsilon_m^k = \frac{\{(\nu)_{m+1}\}^{\beta} p^{m+1}}{\{(1-\varepsilon_m)(m+1)!\}^{\alpha} {}_1S^{\beta}_{\alpha-1,m}(\nu;1;p)}$$

For $\alpha - \beta \geq 1$, this truncated approximation is good because $\varepsilon_m = O(m)$ and thus, the truncation point $m$ is not large. However, for $0 < \alpha - \beta < 1$ and $p > 1$, the truncation point become too large to compute the approximation. For example, when $\nu = 1.5$, $p = 3, \alpha = 3.1$, $\beta = 3$, m has to be over 50000. This is not practicable. To avoid this difficulty it is useful to make a restriction for the parameter $p$ such that $p < 1$ when $\alpha - \beta \to 0$. For example, with the restriction $p < 10^{\alpha - \beta}$, we see the relative truncation error $R_{50}(1.5, 3, 3.1, 3) < 0.001$.

Next using the Laplace's method (Bleistein and Handelesman, 1986, Ch 8.3) we have obtained an asymptotic approximation formula for the normalizing constant ${}_1S^{\beta}_{\alpha-1}(\nu, 1, p)$ as

$$ {}_1\tilde{S}^{\beta-1}_{\alpha-1}(\nu;1;p) = \frac{p^{\{1-\alpha+(2\nu-1)\beta\}/2(\alpha-\beta)} \exp\{(\alpha-\beta) p^{1/(\alpha-\beta)}\}}{(2\pi)^{(\alpha-\beta-1)/2}\sqrt{\alpha-\beta}} \qquad (9)$$

This formula reduces to asymptotic formula by Minka et al. (2003) when $\nu = 1$ or $\beta = 0$ and that by Imoto (2014) when $\alpha = 1$. The formula (9) has been derived for an integer values $\alpha \geq 1$ and $-\beta \geq 0$, but numerical studies suggest that it holds for $0 < \alpha - \beta < 1$ and $p > 1$, where it is difficult to compute ${}_1S^{\beta}_{\alpha-1}(\nu, 1, p)$ by truncated approximation (8).

Taking $m = 18000$ such that $R_m(\nu, p, \alpha, \beta) < 10^{-28}$, Table 1 gives the percentage errors $100\{_1\tilde{S}_{\alpha-1}^{\beta}(\nu;1;p) - {_1}S_{\alpha-1,m}^{\beta}(\nu;1;p)\}/_1 S_{\alpha-1,m}^{\beta}(\nu;1;p)$ .

**Table 1**. The percentage errors for approximation (9) with $\beta = 2.5$

| $p \setminus \alpha$ | $\nu = 0.5$ | | | | | $\nu = 1.5$ | | | | |
|---|---|---|---|---|---|---|---|---|---|---|
| | 2.6 | 2.7 | 2.8 | 2.9 | 3.0 | 2.6 | 2.7 | 2.8 | 2.9 | 3.0 |
| 1.0 | 37 | 2 | -14 | -24 | -31 | -88 | -80 | -73 | -68 | -31 |
| 1.2 | -57 | -44 | -42 | -43 | -44 | -55 | -61 | -60 | -58 | -44 |
| 1.4 | -49 | -58 | -54 | -52 | -51 | -21 | -42 | -47 | -48 | -51 |
| 1.6 | -8 | -55 | -57 | -55 | -54 | -7 | -28 | -37 | -40 | -54 |
| 1.8 | -2 | -38 | -54 | -55 | -55 | -2 | -18 | -28 | -33 | -55 |
| 2.0 | -1 | -18 | -46 | -52 | -54 | -1 | -11 | -22 | -28 | -54 |
| 2.2 | 0 | -9 | -34 | -47 | -51 | 0 | -7 | -17 | -23 | -51 |
| 2.4 | 0 | -5 | -23 | -40 | -47 | 0 | -5 | -13 | -20 | -47 |

Since the ECOMP distribution is a member of exponential family, the mean is given by differentiating the logarithm of the normalizing constant with respect to $p$. Here we consider the differentiation of the logarithm of the function (9), or

$$p^{1/(\alpha-\beta)} + (1 - \alpha + (2\nu - 1)\beta)/(2(\alpha - \beta)) .$$

This function approximates the mean of the ECOMP distribution for large $p$ and small $|\alpha - \beta|$, where it is difficult to compute the approximation by truncation. Table 2 gives the percentage errors of the approximated mean.

**Table 2**. The percentage errors for approximated mean with $\beta = 2.5$

| $p \setminus \alpha$ | $\nu = 0.5$ | | | | | $\nu = 1.5$ | | | | |
|---|---|---|---|---|---|---|---|---|---|---|
| | 2.6 | 2.7 | 2.8 | 2.9 | 3.0 | 2.6 | 2.7 | 2.8 | 2.9 | 3.0 |
| 1.0 | -1055 | -731 | -575 | -475 | -403 | 107 | 78 | 64 | 55 | 48 |
| 1.2 | -187 | -290 | -276 | -249 | -221 | 30 | 39 | 38 | 36 | 34 |
| 1.4 | 46 | -40 | -93 | -107 | -106 | 5 | 18 | 23 | 24 | 24 |
| 1.6 | 1 | 42 | 2 | -24 | -35 | 1 | 8 | 14 | 16 | 17 |
| 1.8 | 0 | 28 | 36 | 19 | 6 | 0 | 3 | 8 | 11 | 13 |
| 2.0 | 0 | 7 | 35 | 35 | 27 | 0 | 1 | 5 | 8 | 10 |
| 2.2 | 0 | 1 | 22 | 35 | 34 | 0 | 1 | 3 | 6 | 7 |
| 2.4 | 0 | 0 | 11 | 29 | 35 | 0 | -0 | 2 | 4 | 6 |

*2.3 Recurrence relation for probabilities*
The ECOMP $(\nu, p, \alpha, \beta)$ pmf has a simple recurrence relation given by

$$\frac{P(X=k+1)}{P(X=k)} = \frac{p(\nu+k)^{\beta}}{(k+1)^{\alpha}} \Rightarrow (k+1)^{\alpha} P(X=k+1) = p(\nu+k)^{\beta} P(X=k) \qquad (10)$$

with $P(X=0) = [_1 S_{\alpha-1}^{\beta}(\nu;1;p)]^{-1}$. This will be useful for the computation of the probabilities.

### 2.4 Exponential family

The pmf in (6) can also be expressed as

$$P(X=k) = \exp[\beta \log(\nu)_k - \beta \log \Gamma(k) - \alpha \log k! + k \log p - \log_1 S_{\alpha-1}^{\beta}(\nu;1;p)] \qquad (11)$$

Which immediately implies that the ECOMP $(\nu, p, \alpha, \beta)$ pmf belongs to the *exponential family* with parameters $(\log p, \alpha, \beta)$ when $\nu$ is a nuisance parameter or its value <span style="color:red">is given.</span>

### 2.5 Dispersion level of with respect to the Poisson distribution

The pmf in (6) can be seen as a weighted Poisson $(p)$ distribution with weight function

$$w(x) = \{\Gamma(\nu+x)\}^{\beta}/(\Gamma(1+x))^{\alpha-1}$$
$$= \exp[\beta \log \Gamma(\nu+x) - (\alpha-1) \log \Gamma(1+x)]$$
$$= \exp[\{\beta \log \Gamma(\nu+x) + (1-\alpha) \log \Gamma(1+x)\}]$$

Now $\log \Gamma(\nu+x)$ is known to be convex function of $x$, which implies $\log \Gamma(1+x)$ is also a convex function of $x$. Therefore if $\alpha < 1$ and $\beta > 0$, $\beta \log \Gamma(\nu+x) + (1-\alpha) \log \Gamma(1+x)$ is also convex function. But it is known from Castillo and Perez-Casany (2005) that weighted Poisson with weight function $w(x) = \exp[c.f(x)]$ where $f(x)$ is a convex function is over (under) dispersed if $c > (<)0$.

Therefore ECOMP $(\nu, p, \alpha, \beta)$ is over dispersed for $1 > \alpha \geq \beta > 0$ that is for $0 < \beta \leq \alpha \leq 1$ and under dispersed for $\{\alpha \leq 1, \beta < 0\} \bigcup \{\alpha > 1, \beta > 0\} \bigcup \{\alpha > 1, \beta < 0\}$.

As particular cases of the above result, we can see that the COM-NB distribution (i.e. when $\beta = 1$) is always over dispersed for $\alpha = 1$ (which is when it reduces to Negative binomial) and under dispersed for $\alpha > 1$. Again GCOMP distribution (i.e. when $\alpha = 1$) is over dispersed for $0 < \beta \leq 1$ and under dispersed for $\beta < 0$. The new generalized NB distribution with pmf (7) is over dispersed when $\gamma = 1$ (which is when it reduces to Negative binomial) and under dispersed if $\gamma > 1$.

### 2.6 Dispersion level of with respect to the COM-Poisson and GCOM-Poisson distribution

Following the same steps as above it can be checked that w.r.t COM-Poisson, ECOMP $(\nu, p, \alpha, \beta)$ is over dispersed for $\alpha \geq \beta > 0$ and under dispersed for $\beta < 0$ and w.r.t GCOM-Poisson it is over dispersed for $\beta \leq \alpha < 1$ and under dispersed for $1 < \beta \leq \alpha$.

### 2.7 Different formulations of ECOMP $(\nu, p, \alpha, \beta)$

Two different formulations of the proposed distribution are presented in this section.

### 2.7.1 ECOMP $(\nu, p, \alpha, \beta)$ as a distribution from a queuing set up

Like the COM-Poisson distribution (Conway and Maxwell, 1962), ECOMP $(\nu, p, \alpha)$ can also be derived as the probability of the system being in the $k^{\text{th}}$ state for a queuing system with state dependent service and arrival rate.

Consider a single server queuing system with state dependent (that is dependent on the system state, $k^{th}$ state means $k$ number of units in the system) arrival rate $\lambda_k = (\nu + k)^\beta \lambda, k \geq 1$, and state dependent service rate $\mu_k = k^\alpha \mu, k \geq 1$, where, $1/\mu$ and $1/\lambda$ are respectively the normal mean service and mean arrival time for a unit when that unit is the only one in the system; $\alpha$ and $\nu$ are the pressure coefficients, reflecting the degree to which the service and arrival rates of the system are affected by the system state. This set up implies that while the arrival rate and the service rate increases exponentially as queue lengthens (i.e. $n$ increases).

Following Conway and Maxwell (1962), the system differential difference equations are given by

$$P_0(t+\Delta) = (1 - \lambda \nu^\beta \Delta) P_0(t) + \mu \Delta P_1(t) \qquad (12)$$

and $P_k(t+\Delta) = (1 - \lambda(\nu+k)^\beta \Delta - \mu k^\alpha \Delta) P_k(t) + \lambda(\nu+k-1)^\beta \Delta P_{k-1}(t)$

$$+ \mu(k+1)^\alpha \Delta P_{k+1}(t)) \qquad k = 1, 2, \cdots \qquad (13)$$

Let $\lambda/\mu = p$. Then from (12) and (13) we get

$\{P_0(t+\Delta) - P_0(t)\}/\Delta = -\mu p \nu^\beta P_0(t) + \mu P_1(t)$ and

$\{P_k(t+\Delta) - P_k(t)\}/\Delta = \mu p(\nu+k-1)^\beta P_{k-1}(t) - \mu(p(\nu+k) + k^\alpha) P_k(t)$

$+ \mu(k+1)^\alpha P_{k+1}(t)), k = 1, 2, \cdots$

Now as $\Delta \to 0$ we get

$P_0^{/}(t)\} = -\mu p \nu^\beta P_0(t) + \mu P_1(t)$ and

$P_k^{/}(t) = \mu p(\nu+k-1)^\beta P_{k-1}(t) - \mu(p(\nu+k)^\beta + k^\alpha) P_k(t) + \mu(k+1)^\alpha P_{k+1}(t)), k = 1, 2, \cdots$

Assuming a steady state (i.e. $P_k^{/}(t) = 0$ for all $k$) we get

$P_1(t) = p \nu^\beta P_0(t)$ and

$P_{k+1}(t) = \dfrac{k^\alpha}{(k+1)^\alpha} P_k(t) - \dfrac{(\nu+k)^\beta}{(k+1)^\alpha} p P_k(t) + \dfrac{(\nu+k-1)^\beta}{(k+1)^\alpha} p P_{k+1}(t), k = 1, 2, \cdots$

Putting $k = 1$ we get

$$P_2(t) = \dfrac{1}{2^\alpha} P_1(t) - \dfrac{(\nu+1)^\beta}{2^\alpha} p P_1(t) + \dfrac{\nu^\beta z}{2^\alpha} p P_0(t)$$

$$= \dfrac{1}{2^\alpha} \nu^\beta p P_0(t) - \dfrac{(\nu+1)^\beta}{2^\alpha} \nu^\beta p P_0(t) + \dfrac{\nu^\beta z}{2^\alpha} p P_0(t)$$

$$= \dfrac{\{\nu(\nu+1)\}^\beta}{(2!)^\alpha} p^2 P_0(t) = \dfrac{\{(\nu)_2\}^\beta}{(2!)^\alpha} p^2 P_0(t)$$

Similarly, for $k = 2$ we get

$$P_3(t) = \left( \dfrac{2^\alpha}{3^\alpha} - \dfrac{(\nu+2)^\beta}{3^\alpha} p \right) P_2(t) + \dfrac{(\nu+1)^\beta p}{3^\alpha} P_1(t)$$

$$= \left( \frac{2^\alpha}{3^\alpha} - \frac{(\nu+2)^\beta}{3^\alpha} p \right) \frac{\{\nu(\nu+1)\}^\beta}{(2!)^\alpha} p^2 P_0(t) + \frac{(\nu+1)^\beta z}{3^\alpha} \nu \, p \, P_0(t)$$

$$= \frac{\{\nu(\nu+1)(\nu+2)\}^\beta}{(3!)^\alpha} p^3 P_0(t)$$

In general,

$$P_k(t) = \frac{\{(\nu)_k\}^\beta}{(k!)^\alpha} p^k P_0(t), \text{ where } P_0(t) = 1 / \sum_{i=0}^\infty \left\{ \frac{\{(\nu)_i\}^\beta}{(i!)^\alpha} p^i \right\}.$$

Since we have assumed a steady state (i.e. $P_k^/(t) = 0$ for all $k$) $P_k(t)$ can be replaced by $P_k$.

### 2.7.2 ECOMP $(\nu, p, \alpha, \beta)$ as exponential combination formulation

The general form of the exponential combination of two pmfs say $f_1(x; \theta_1)$ and $f_2(x; \theta_2)$ is given by (Atkinson, 1970)

$$\{f_1(x; \theta_1)\}^\beta \{f_2(x; \theta_2)\}^{1-\beta} / \sum f_1(x; \theta_1)^\beta f_2(x; \theta_2)^{1-\beta}$$

This combining of the pmf was suggested by Cox (1961, 1962) for combining the two hypotheses ($\beta = 1$ i.e. the distribution is $f_1$ and $\beta = 0$ that is the distribution is $f_2$) in a general model of which they would both be special cases. The inferences about $\beta$ made in the usual way and testing the hypo- thesis that the value of $\beta$ is zero or one is equivalent to testing for departures from one model in the direction of the other.

Now the probability function resulting from the exponential combination of NB $(\nu, \lambda)$ and COM-Poisson $(\mu, \theta)$ is given by

$$\left\{ \frac{(\nu)_k}{k!} \lambda^k \right\}^\beta \left\{ \frac{\mu^k}{(k!)^\theta} \right\}^{1-\beta} \Bigg/ \sum_{i \geq 0} \left\{ \frac{(\nu)_i}{i!} \lambda^i \right\}^\beta \left\{ \frac{\mu^i}{(i!)^\theta} \right\}^{1-\beta}$$

$$= \frac{\{(\nu)_k\}^\beta \{\lambda^\beta \mu^{1-\beta}\}^k}{(k!)^{\theta(1-\beta)+\beta}} \Bigg/ \sum \frac{\{(\nu)_i\}^\beta \{\lambda^c \mu^{1-\beta}\}^i}{(i!)^{\theta(1-\beta)+\beta}}$$

$$= \frac{\{(\nu)_k\}^\beta p^k}{(k!)^\alpha} \Bigg/ \sum_{i \geq 0} \frac{\{(\nu)_i\}^\beta p^i}{(i!)^\alpha}, \text{ substituting } \lambda^\beta \mu^{1-\beta} = p \text{ and } \alpha = \theta(1-\beta) + \beta$$

This is the pmf of ECOMP $(\nu, p, \alpha, \beta)$.

From the above formulations it is clear that ECOMP $(\nu, p, \alpha, \beta)$ can be regarded as a natural extension of COM-Poisson of Conway and Maxwell (1962), COM-NB of Chakraborty and Ong (2014), as well as the GCOMP distribution of Imoto (2014).

### 2.8 Log-concavity and modality

➢ The ECOMP $(\nu, p, \alpha, \beta)$ has *a log-concave probability mass function* when $\{\nu \geq 1, p > 0, \alpha > \beta\}$ since for this distribution (Gupta et al., 1997)

$$\Delta \eta(t) = \frac{P(t+1)}{P(t)} - \frac{P(t+2)}{P(t+1)} = p \frac{(\nu+t)^\beta (t+2)^\alpha - (\nu+t+1)^\beta (t+1)^\alpha}{(t+1)^\alpha (t+2)^\alpha} > 0, \text{ for all } t.$$

From this result the corresponding results of COM-NB $(\nu, p, \alpha)$ and GCOMP $(\nu, p, \beta)$ can be obtained as particular cases. That is COM-NB $(\nu, p, \alpha)$ is log-concave when $\{\nu \geq 1, p > 0, \alpha > 1\}$ and GCOMP $(\nu, p, \beta)$ is log-concave when $\{\nu \geq 1, p > 0, \beta < 1\}$.

Following two important results follows as a consequence of log-concavity:

If $\{\nu \geq 1, p > 0, \alpha > \beta\}$ the ECOMP $(\nu, p, \alpha, \beta)$ distribution is

- a *strongly unimodal* distribution
- has an *increasing failure rate* function

➢ Using the recurrence relation of the probabilities in (10) it can be shown easily that the ECOMP $(\nu, p, \alpha, \beta)$ has

(i) *a non increasing pmf with a unique mode at $X = 0$ if $\nu^{\beta} p < 1$*,
   *e.g.* $\nu = 2, \alpha = 3, \beta = 2$, $p$ should be less than 0.25 to have unique mode at $X = 0$.

(ii) *a unique mode at $X = k$ if $k^{\alpha} / (\nu + k - 1)^{\beta} < p < (k + 1)^{\alpha} / (\nu + k)^{\beta}$*
   *e.g.* $\nu = 2, \alpha = 3, \beta = 2$, $p$ should be between 1.6875 and 2.560 to have unique mode at $X = 3$.

(iii) *two modes at $X = k$ and $X = k - 1$ if $(\nu + k - 1)^{\beta} p = k^{\alpha}$*. In particular the two modes are at $X = 0$ and $X = 1$ if $\nu^{\beta} p = 1$.
   *e.g.* $\nu = 2, \alpha = 3, \beta = 2$, $p$ should be equal to 4.408 to have two modes at $X = 5$ and $X = 6$.

Graphical illustrations of the above three examples are presented in figure 1 by plotting the corresponding pmfs.

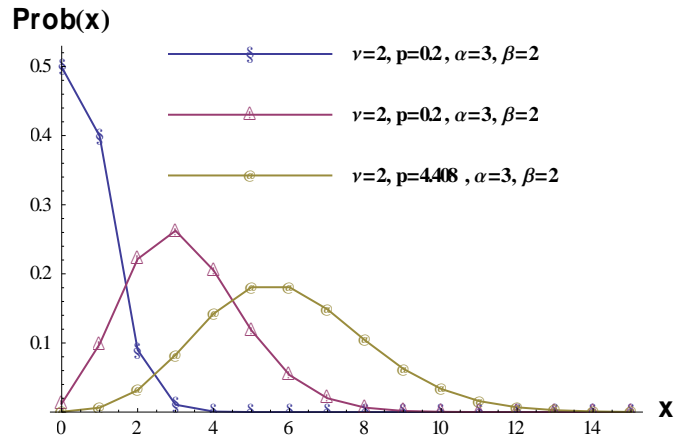

**Figure 1.** pmfs of ECOMP

**Remark3**. For $0 < \nu < 1$ the distribution may be bimodal with one of the modes always at zero. A few instances are presented in figure 2.

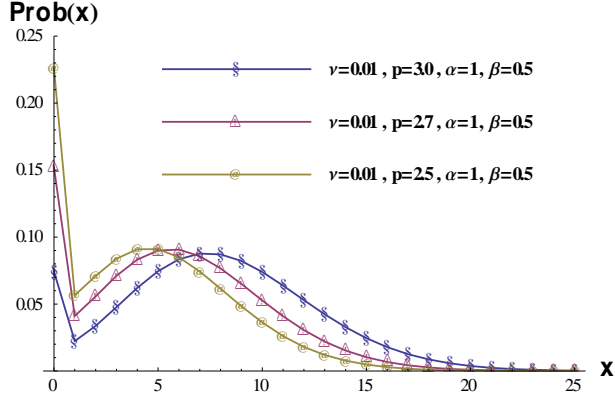

**Figure 2**. Bimodal pmfs of ECOMP

➤ ECOMP $(\nu, p, \alpha, \beta)$ has *a log-convex probability mass function* for $\{0 < \nu \le 1, \alpha = \beta\}$

Proof: ECOMP $(\nu, p, \alpha, \beta)$ has *a log-convex probability mass function* if $\Delta \eta(t) \le 0$, for all $t$. That is if

$$(\nu + t)^{\beta} (t + 2)^{\alpha} < (\nu + t + 1)^{\beta} (t + 1)^{\alpha}$$

$$\Rightarrow \{(t + 2)/(t + 1)\}^{\alpha} < \{(\nu + t + 1)/(\nu + t)\}^{\beta}$$

$$\Rightarrow \{(1 + 1/(t + 1)\}^{\alpha} < \{1 + 1/(\nu + t)\}^{\beta} \tag{14}$$

Since $\alpha \ge \beta$, the condition (14) can hold only for $\nu \le 1$.

Again for $0 < \nu \le 1$ the condition (14) implies

$$\Rightarrow \alpha / \beta \le \log(1 + 1/(t + \nu))/\log(1 + 1/(t + 1)), \text{ for all } t$$

But minimum of $\{ \log(1 + 1/(t + \nu))/\log(1 + 1/(t + 1)) \}$ w.r.t. $t$ is 1

$\Rightarrow \alpha / \beta \le 1$. This implies $\alpha = \beta$ since $\alpha \ge \beta$

*2.9 Momemts*

The $r^{\text{th}}$ factorial moment $E(X^{[r]}) = \mu^{[r]}$ of the ECOMP $(\nu, p, \alpha, \beta)$ is given by

$$\mu^{[r]} = \frac{\{(\nu)_r\}^{\beta} p^r}{(r!)^{\alpha - 1}} \frac{{}_1 S_{\alpha-1}^{\beta}(\nu + r, r + 1, p)}{{}_1 S_{\alpha-1}^{\beta}(\nu, 1, p)} = \frac{\{(\nu)_r\}^{\beta} p^r}{(r!)^{\alpha - 1}} \frac{{}_\beta F_{\alpha-1}(\nu + r; r + 1, r + 1, \cdots, r + 1; p)}{{}_\beta F_{\alpha-1}(\nu; 1, 1, \cdots, 1; p)},$$

# 3. Reliability characteristics and stochastic ordering

*3.1 Survival and Failure rate functions*

The survival function is given by

$$S(t) = 1 - P(X < t) = 1 - \frac{1}{{}_1 S_{\alpha-1}^{\beta}(\nu, ; 1; p)} \sum_{k=0}^{t-1} \frac{\{(\nu)_k\}^{\beta} p^k}{(k!)^{\alpha}}$$

$$= 1 - \frac{1}{{}_\beta F_{\alpha-1}(\nu; 1, 1, \cdots, 1; p)} \sum_{k=0}^{t-1} \frac{\{(\nu)_k\}^{\beta} p^k}{(k!)^{\alpha}}$$

$$S(t) = P(X \geq t) = \frac{(\nu)_t \ p^t}{(t!)^\alpha} \cdot \frac{{}_{\beta+1}F_{\alpha-1}(\nu+t,1; \ t+1,t+1,\cdots,t+1; \ p)}{{}_\beta F_{\alpha-1}(\nu; 1,1,\cdots,1; \ p)},$$

while the failure rate function is

$$r(t) = \frac{P(X = t)}{P(X \geq t)} = \frac{1}{{}_2 S_{\alpha-1}^\beta(\nu+t,1; \ t+1; \ p)} = \frac{1}{{}_{\beta+1}F_\alpha(\nu+t,1; \ t+1,t+1,\cdots,t+1; \ p)},$$

where the second expressions in terms of hypergeometric function is for the case when $\alpha, \beta$ are positive integers.

## 3.2 Stochastic orderings

**Theorem 1.** The ECOMP $(\nu, p, \alpha, \beta)$ is smaller than the COM-NB $(\nu, p, \alpha)$ distribution in the likelihood ratio order i.e. $X \leq_{lr} Y$ when $\beta > 1$.

**Proof:** If $X \sim$ ECOMP$(\nu, p, \alpha, \beta)$ and $Y \sim$ COM-NB$(\nu, p, \alpha)$, then

$$\frac{P(Y = n)}{P(X = n)} = \{(\nu)_n\}^{1-\beta} \frac{{}_1 S_{\alpha-1}^\beta(\nu,1,p)}{{}_1 S_{\alpha-1}^1(\nu,1,p)}.$$

This is clearly increasing in $n$ as $\beta > 1$ (Shaked and Shanthikumar, 2007 and Gupta et al., 2014). Hence the result is proved.

**Theorem 2.** The ECOMP $(\nu, p, \alpha, \beta)$ is smaller than the GCOMP $(\nu, p, \beta)$ distribution in the likelihood ratio order i.e. $X \leq_{lr} Y$ when $\alpha > 1$.

**Proof:** If $X \sim$ ECOM-NB$(\nu, p, \alpha, \beta)$ and $Y \sim$ GCOMP$(\nu, p, \beta)$, then

$$\frac{P(Y = n)}{P(X = n)} = \frac{(n!)^{\alpha-1}}{(\Gamma \nu)^\beta} \frac{{}_1 S_{\alpha-1}^\beta(\nu,1,p)}{{}_1 S_0^\beta(\nu,1,p)}.$$

This is clearly increasing in $n$ as $\alpha > 1$ (Shaked and Shanthikumar, 2007 and Gupta et al., 2014). Hence the result is proved.

## 4. A numerical example

To fit the proposed distribution, we have to estimate the parameters $(\nu, p, \alpha, \beta)$ in (6). The maximum likelihood (ML) estimation is often used for fitting to real data, but the log likelihood function of the proposed distribution

$$L(\nu, p, \alpha, \beta) = \beta \sum_{i=1}^k f_i \log(\nu)_i - \alpha \sum_{i=1}^k f_i \log i! + \log p \sum_{i=1}^k i \ f_i - N \log_1 S_{\alpha-1}^\beta(\nu,1,p) \qquad (15)$$

where $f_i$ is the observed frequency of $i$ events, $N = \sum_{i=1}^k f_i$, $k$ is the highest observed value,

has some local maximum points for some datasets, or the likelihood equations do not always have unique solution.

Here we have used profile likelihood method to fit the proposed distribution to the data of Corbet (1942) on Malayan butteries with zeros, which also has been used by Blumer (1974) for fitting the Poisson-lognormal distribution

$$P(x) = \frac{1}{x!} \int_0^\infty \frac{e^{-z} z^x}{z\sqrt{2\pi\sigma^2}} \exp\{-(\log z - \nu)^2 / 2\sigma^2\} dz$$

Corbet caught altogether 620 species, but he also estimated that the total buttery fauna of the area contained 924 species, so that 304 species were missing from the collection and treated as count zero. We also fit the Poisson-lognormal (P-Log), COMNB ((2) of section 1) and GCOMP ((3) of Section 1) distributions to the same data for comparing the fitting with the proposed distribution. The ML estimates and fitting results of the distributions are given in Table 3.

**Table 3**: Distribution of Corbet's Malayan Buttery with zeros (Corbet, 1942)

| Count | Observed | P-Log | COMNB | GCOMP | ECOMP |
|---|---|---|---|---|---|
| 0 | 304 | 294.75 | 315.36 | 315.36 | 304.97 |
| 1 | 118 | 127.35 | 94.24 | 94.24 | 117.12 |
| 2 | 74 | 74.58 | 59.76 | 59.76 | 67.25 |
| 3 | 44 | 50.73 | 44.58 | 44.58 | 45.92 |
| 4 | 24 | 37.54 | 35.74 | 35.74 | 34.51 |
| 5 | 29 | 29.29 | 29.85 | 29.85 | 27.57 |
| 6 | 22 | 23.71 | 25.60 | 25.60 | 22.94 |
| 7 | 20 | 19.71 | 22.37 | 22.37 | 19.66 |
| 8 | 19 | 16.73 | 19.81 | 19.81 | 17.23 |
| 9 | 20 | 14.43 | 17.73 | 17.73 | 15.35 |
| 10 | 15 | 12.61 | 16.00 | 16.00 | 13.86 |
| 11 | 12 | 11.14 | 14.53 | 14.53 | 12.64 |
| 12 | 14 | 9.93 | 13.27 | 13.27 | 11.63 |
| 13 | 6 | 8.93 | 12.18 | 12.18 | 10.77 |
| 14 | 12 | 8.08 | 11.23 | 11.23 | 10.04 |
| 15 | 6 | 7.35 | 10.38 | 10.38 | 9.39 |
| 16 | 9 | 6.73 | 9.63 | 9.63 | 8.83 |
| 17 | 9 | 6.19 | 8.96 | 8.96 | 8.32 |
| 18 | 6 | 5.71 | 8.35 | 8.35 | 7.86 |
| 19 | 10 | 5.29 | 7.80 | 7.80 | 7.45 |
| 20 | 10 | 4.92 | 7.30 | 7.30 | 7.07 |
| 21 | 11 | 4.59 | 6.84 | 6.84 | 6.73 |
| 22 | 5 | 4.29 | 6.42 | 6.42 | 6.40 |
| 23 | 3 | 4.02 | 6.04 | 6.04 | 6.10 |
| 24 | 3 | 3.78 | 5.69 | 5.69 | 5.82 |
| 25+ | 119 | 131.64 | 114.36 | 114.36 | 118.57 |
| Total | 924 | 924 | 924 | 924 | 924 |
| AIC | | 4518.88 | 4518.22 | 4518.22 | **4510.02** |
| $\chi^2$ | | 36.86 | 28.50 | 28.50 | **18.57** |
| $p$-value | | 0.03 | 0.15 | 0.15 | **0.61** |
| MLEs | | $\hat{\mu}=0.70$ | $\hat{v}=0.31$ | $\hat{v}=0.31$ | $\hat{v}=2.86$ |
| | | $\hat{\sigma}^2=5.42$ | $\hat{p}=0.97$ | $\hat{p}=0.97$ | $\hat{p}=1.26$ |
| | | | $\hat{\alpha}=1.00$ | $\hat{\beta}=1.00$ | $\hat{\alpha}=-1.07$ |
| | | | | | $\hat{\beta}=-1.13$ |

For this dataset, the ML estimate $\hat{\alpha}$ of the COMPNB distribution and ML estimate $\hat{\beta}$ of the GCOMP distribution show these two distributions reduce to the negative binomial distribution while the proposed distribution does not seem to reduce to the negative binomial distribution. Actually, the likelihood ratio test for $H_0$: Negative binomial distribution ($\alpha = \beta = 1$) Vs $H_1$: ECOMP distribution ($\alpha \neq 1$ or $\beta \neq 1$) rejects the negative binomial distribution ($p$-value is 0.001).

The Poisson-lognormal distribution was derived for long-tailed count data and fitted to the same data as in this section. Comparing the distributions, we see that the ECOMP distribution gives better fitting in the sense of AIC and $\chi^2$ goodness of fit, and especially gives very good fittings for the count 0 and the tail part 25+. We have tried fitting the ECOMP distribution to various count data and seen that the distribution gives good fittings for the count data with many zeros or with long-tail.

## 5. Concluding Remarks

In this paper a new discrete distribution that extends the Conway-Maxwell-Poisson distribution and also unifies the COM-NB (Chakraborty and Ong, 2014) and GCOMP (Imoto, 2014) is proposed and its important distributional properties investigated. This distribution which arises from queuing theory set up and also as exponential combination has many desirable properties with potential applications in modeling varieties of count data.

## References

1. Atkinson, A. C. (1970) A Method for Discriminating Between Models. Journal of the Royal Statistical Society. Series B (Methodological), **32**, 3, pp. 323-353

2. Blumer, M. G. (1974). On the fitting the Poisson-lognormal distribution to species abundance data. Biometrics: 30, 101-110.

3. Bleistein, N., Handelsman, R. (1986) *Asymptotic Expansions of Integrals*, Dover, NY.

4. Castillo, J., Pérez-Casany, M. (2001) Overdispersed and underdispersed Poisson generalizations, J. Stat. Plan. Inference, **134** (2005) 486–500.

5. Chakraborty, S. and Ong, S. H. (2014), A COM-type Generalization of the Negative Binomial Distribution, Accepted in April 2014, to appear in Communications in Statistics-Theory and Methods

6. Conway, R. W., Maxwell W. L. (1962) A queueing model with state dependent service rates. J Industrl Engng **12**:132–136

7. Corbet, A. S. (1942). The distribution of butteries in the Malay Peninsula. Proceedings of the Royal Entomological Society of London, Series A, General Entomology: **16**, 101-116.

8. Cox, D. R. (1961). Tests of separate families of hypotheses. Proc. 4th Berkeley Symp., **1**, 105-123.

9. Cox, D. R. (1962). Further results on tests of separate families of hypotheses. J. R. Statist. Soc. B, **24**: 406-424.

10. Gupta, P. L., Gupta, R. C., Tripathi, R. C. (1997) On the monotonic properties of discrete failure rates. Journal of Statistical Planning and Inference 65:255–268

11. Gupta R. C., Sim S. Z., Ong S. H. (2014) Analysis of Discrete Data by Conway-Maxwell Poisson Distribution. AStA Advances in Statistical Analysis, DOI 10.1007/s10182-014-0226-4

12. Imoto, T. (2014) A generalized Conway-Maxwell-Poisson distribution which includes the negative binomial distribution, Applied Mathematics and Computation, **247**, 824-834

13. Johnson N. L., Kemp A. W., Kotz, S. (2005) *Univariate Discrete Distributions*, Wiley, New York

14. Minka T. P., Shmueli G., Kadane J. B., Borle S., and Boatwright P. (2003). Computing with the COM-Poisson distribution. Technical Report: 776, Department of Statistics, Carnegie Mellon University, http://www.stat.cmu.edu/tr/tr776/tr776.html.

15. Shaked M., Shanthikumar J. G. (2007*)* Stochastic orders. Springer Verlag